%% file: kronsum-arxiv.tex
\documentclass[11pt]{article}
\RequirePackage[OT1]{fontenc}
\RequirePackage{amsthm,amsmath}
\RequirePackage{hyperref}
\usepackage{amsfonts}
\usepackage{amssymb}
\usepackage{amstext}
\usepackage{parskip}
\usepackage{fullpage}
\usepackage{graphicx}
\usepackage{subcaption}

\RequirePackage[numbers]{natbib}
\usepackage{amsfonts,amssymb,amsthm,amsmath,enumerate}

\newcommand{\be}{\begin{eqnarray*}}
\newcommand{\ee}{\end{eqnarray*}}
\newcommand{\ben}{\begin{eqnarray}}
\newcommand{\een}{\end{eqnarray}}

\newcommand{\var}{{\rm Var}}
\newcommand{\cov}{{\rm Cov}}
\newcommand{\tr}{{\rm tr}}
\newcommand{\argmin}{\operatornamewithlimits{argmin}}

\newcommand{\grad}{\nabla}

\newcommand{\E}{\mathbb{E}}

\newcommand{\twonorm}[1]{\left\lVert#1\right\rVert_2}

\newcommand{\shtwonorm}[1]{\lVert#1\rVert_2}
\newcommand{\shonenorm}[1]{\lVert#1\rVert_1}
\newcommand{\fnorm}[1]{\left\lVert#1\right\rVert_F}
\newcommand{\norm}[1]{\left\lVert#1\right\rVert}
\newcommand{\abs}[1]{\left\lvert#1\right\rvert}

\newcommand{\minus}{\ensuremath{-}}

\newcommand{\M}{{\mathcal{M}}}

\newcommand{\V}{\mathcal{V}}

\newcommand{\N}{{\mathcal N}}

\newcommand{\ora}{{\rm oracle}}

\newcommand{\inv}[1]{\frac{1}{#1}}
\def\conv{\mathop{\text{\rm conv}\kern.2ex}}

\newcommand{\RE}{\textnormal{\textsf{RE}}}
\newcommand{\ip}[1]{\;\langle{\,#1\,}\rangle\;}

\newcommand{\maxnorm}[1]{\ensuremath{\left|#1\right|_{\max}}}
\newcommand{\expct}[1]{\ensuremath{\mathbb E}\left[#1\right]}
\newcommand{\silent}[1]{}
\newcommand{\mvec}[1]{\rm{vec}\left\{\,#1\,\right\}}

\newcommand{\ve}{\varepsilon}

\def\qed{\hskip1pt $\;\;\scriptstyle\Box$}
\def\Ber{\mathop{\text{Bernoulli}\kern.2ex}}
\def\supp{\mathop{\text{supp}\kern.2ex}}

\def\corr{\mathop{\text{corr}\kern.2ex}}
\def\prec{\mathop{\text{precision}\kern.2ex}}
\def\recall{\mathop{\text{recall}\kern.2ex}}
\def\mnorm{\mathcal{N}_{f,m}\kern.2ex}
\def\var{\mathop{\text{Var}\kern.2ex}}
\def\ess{\mathop{\text{ess}\kern.2ex}}
\def\dom{\mathop{\text{dom}\kern.2ex}}
\def\lin{\mathop{\text{lin}\kern.2ex}}

\newcommand{\func}[1]{\ensuremath{\mathrm{#1}}}
\newcommand{\diag}{\func{diag}}

\newenvironment{proofof}[1]{\vskip5pt \hspace*{20pt}{\it Proof}{ of
    #1}.\hskip10pt}{\qed\vskip5pt}

{\end{eqnarray}\end{subequations}\hskip-4.0pt}

\let\hat\widehat
\let\tilde\widetilde

\def\supp{\mathop{\text{\rm supp}\kern.2ex}}
\def\argmin{\mathop{\text{arg\,min}\kern.2ex}}

\newcommand{\prob}[1]{\ensuremath{\mathbb P}\left(#1\right)}

\newcommand{\beq}{\begin{equation}}
\newcommand{\eeq}{\end{equation}}
\newcommand{\bnum}{\begin{enumerate}}
\newcommand{\enum}{\end{enumerate}}
\newcommand{\bit}{\begin{itemize}}
\newcommand{\eit}{\end{itemize}}
\newcommand{\bens}{\begin{eqnarray*}}
\newcommand{\eens}{\end{eqnarray*}}

\newcommand{\R}{{\mathbb R}}

\newcommand{\B}{\ensuremath{\mathcal B}}

\newcommand{\SM}{\ensuremath{{\mathcal{S}^m}}}

\newcommand{\e}{\epsilon}

\newcommand{\vp}{\varphi}

\newcommand{\onem}{\textstyle \frac{1}{m}}
\newcommand{\onen}{\textstyle \frac{1}{n}}

\newcommand{\ignore}[1]{}{}

\newcommand*{\ttilde}{{\raise.17ex\hbox{$\scriptstyle\sim$}}}

\makeatletter
\newsavebox{\mybox}\newsavebox{\mysim}
\newcommand*{\distas}[1]{%
  \savebox{\mybox}{\hbox{\kern3pt$\scriptstyle#1$\kern3pt}}%
  \savebox{\mysim}{\hbox{$\sim$}}%
  \mathbin{\overset{#1}{\kern\z@\resizebox{\wd\mybox}{\ht\mysim}{$\sim$}}}%
}
\makeatother

\makeatletter
\def\moverlay{\mathpalette\mov@rlay}
\def\mov@rlay#1#2{\leavevmode\vtop{%
   \baselineskip\z@skip \lineskiplimit-\maxdimen
   \ialign{\hfil$\m@th#1##$\hfil\cr#2\crcr}}}
\newcommand*{\charfusion}[3][\mathord]{
  #1{\ifx#1\mathop\vphantom{#2}\fi\mathpalette\mov@rlay{#2\cr#3}}
  \ifx#1\mathop\expandafter\displaylimits\fi}
\makeatother

\usepackage{times}
\usepackage{mdframed}
\usepackage{float}

\newtheorem{theorem}{Theorem}[section]

\newtheorem{lemma}[theorem]{Lemma}
\newtheorem{proposition}[theorem]{Proposition}

\newtheorem{definition}[theorem]{Definition}

\newtheorem{remark}[theorem]{Remark}

\def\qed{\hskip1pt $\;\;\scriptstyle\Box$}

\begin{document}

\title{Kronecker sum covariance models for spatio-temporal data}

\author{Shuheng Zhou  \ \ \ \ \ \ \ \ \ \ \ \ \ \ \ \ \ \ \ Seyoung Park \ \ \  \ \ \ \ \ \ \ \ \ \ \ \ \ \ \  Kerby Shedden \\
    University of California, Riverside  \ \ \ \ \  Yonsei University  \ \ \ \ \  University of Michigan, Ann Arbor}

\date{}

\maketitle

\begin{abstract}
In this paper, we study the subgaussian matrix variate model, where we
observe the matrix variate data $X$ which consists of a signal matrix $X_0$ and
a noise matrix $W$. More specifically, we study a subgaussian model using the Kronecker sum
covariance as in Rudelson and Zhou (2017).
Let $Z_1, Z_2$ be independent copies of a subgaussian random matrix  
$Z =(Z_{ij})$ (cf.~\eqref{eq::Zdef}), where $Z_{ij}, \forall i, j$ are independent  
mean 0, unit variance, subgaussian random variables with bounded
$\psi_2$ norm.
We use $X \sim \M_{n,m}(0,  A \oplus B)$
to denote the subgaussian random matrix $X_{n \times m}$ which is 
generated using \eqref{eq::addmodel},
\ben
\label{eq::addmodel}
X  & = & Z_1 A^{1/2} + B^{1/2} Z_2.
\een
In this covariance model, the first component $A \otimes I_n$
describes the covariance of the {\it signal} $X_0 = Z_1 A^{1/2}$,
which is an ${n \times m}$ random design matrix with independent
subgaussian row vectors, and the other component $I_m \otimes B$
describes the covariance for the {\it noise matrix } $W =B^{1/2} Z_2$,
which contains independent subgaussian column vectors $w^1, \ldots,
w^m$, independent of $X_0$.
This leads to a non-separable class of 
models for the observation $X$, which we denote by
$X \sim \M_{n,m}(0, A \oplus B)$ throughout this paper.
Our method on inverse covariance estimation corresponds to the
proposal in~\citet{Yuan10} and \citet{LW12}, only now dropping the
i.i.d. or Gaussian assumptions.
We present the statistical rates of convergence.
\end{abstract}

\input{main}

\section*{Acknowledgment}
Part of this work was conducted while Shuheng Zhou was visiting the Simons Institute
for the Theory of Computing.

\bibliography{subgaussian2}

\appendix

\input{append}

\end{document}

%% file: main.tex
\section{Introduction}
Graphical models capture conditional independence relationships among
variables and are well established under the
i.i.d. assumption. However, graphical models for matrix-variate data
remain relatively underdeveloped. The need for such an extension
arises in scientific experiments, where data corruption and
measurement errors cannot be ignored.
For example, neuroscience experiments often involve a large number of
trials across varying experimental conditions on a relatively small
number of subjects. As another example, motor control in animals is a
complex process requiring the coordination of multiple muscles to
produce intricate movements in space.
Analyzing neurocontrol of flight, therefore, requires the study of high-dimensional temporal 
measurements, which are often significantly affected by measurement
errors, and / or spatial correlations.

In this paper, we introduce a novel Kronecker 
sum covariance structure to model the complex dependencies among rows and 
columns for the subgaussian matrix-variate data $X$ with measurement errors.
In particular, we model such $n\times m$ data $X$ with corrupted entries as 
\begin{equation}
\label{eq:model_intro0}
    X = X_0 + W \sim \mathcal{N}_{n\times m}(0, A \oplus B),
\end{equation}
Here $X_0$ is an $n\times m$ design matrix with dependent columns and
independent rows, and $W$ is a mean-zero $n \times m$ random noise
matrix with independent columns and dependent rows. We assume $X_0$
and $W$ are independent from each other.
Here, $\oplus$ denotes Kronecker sum, i.e. $A\oplus B = A\otimes I_n + I_m
\otimes B$, where $\otimes$ denotes the Kronecker product.
We use $X\sim \mathcal{N}_{n\times m}(0, A\oplus B)$ to denote $\mvec{X} \sim
N(0, A\oplus B)$. Here $\mvec{X}$ is formed by stacking the columns of $X$ into a vector
in $\R^{mn}$. Let $(a)_+ := a \vee 0$ and the Frobenius norm be  $\norm{X}^2_F = \sum_i\sum_j x_{ij}^2$. 

The key difference between our framework and 
the existing work is that we assume that only one copy of
the random matrix $X$ as in~\eqref{eq:model_intro0} is observed (with the single measurement error matrix $W$).
Unlike the single-sample Kronecker product covariance models studied in \citet{Zhou14a} and
\citet{HSZ15}, where the design matrix can be decorrelated from the estimated covariance matrices, we cannot
directly eliminate the noise matrix $W$ from the observed data $X$.
However, statistical methods that utilize the second-order
statistics of $X$, such as regression models and graph estimation, can
still benefit from the Kronecker sum model. Our approach leverages the
column independence of the noise matrix $W$ to effectively mitigate
the impact of noise in the following sense.
For the Kronecker sum model~\eqref{eq::addmodel} to be
identifiable, we assume the trace of $A$ is known as in \citet{RZ15}.
We now state $(A1)$ followed by an initial estimator of $A$:
\bnum
\item[(A1)]
  We assume $\tr(A) = m$ is a known parameter, where $\tr(A)$ denotes
the trace of matrix $A$. By fixing $\tr(A)$, we can construct an estimator for $\tr(B)$:
\ben
\label{eq::trBest}
&& \hat\tr(B) :=
\onem (\fnorm{X}^2 -n \tr(A))_{+},  \hat\tau_B  := \onen \hat\tr(B);\\
\label{eq::hatGamma}
&& \text{Set the Gram matrix} \; \hat\Gamma := \onen X^T X -\hat\tau_B I_{m}.
\een
\enum
As a quick remark, any $m \times m$ correlation matrix satisfies (A1).
Assuming $\tr(A)$ or $\tr(B)$ is known is unavoidable as the covariance model is not identifiable 
otherwise for the single sample case as considered in the present
work. We further develop multiple regression methods using the corrected
Lasso to estimate the inverse of $A$ and $B$ and show statistical rate of convergence.
We will estimate the inverse covariances $\Theta = A^{-1}$ and $\Phi=B^{-1}$ using 
the nodewise regression procedure \citep{MB06} followed by the refit
procedures, only now dropping the i.i.d. or Gaussian assumptions.

\noindent{\bf  The additive Gaussian graphical models.}
To ease discussion,  we first consider~\eqref{eq::addmodel},
where   $Z_1$ and $Z_2$ are assumed to be independent
copies of a Gaussian random ensemble $Z$ with independent standard
normal entries.
Recall when $Z_1$ in~\eqref{eq::addmodel} is a Gaussian random
ensemble, with i.i.d. $N(0, 1)$ entries,  we say the random matrix $X_0$ as
defined in~\eqref{eq::addmodel} follows the matrix-variate normal distribution
\ben 
\label{eq::matrix-normal-rep-intro}
X_{0, n \times m} \sim \N_{n,m}(0, A_{m \times m} \otimes I_{n}). 
\een
This is equivalent to say $\mvec{X_0}$ follows a multivariate normal
distribution with mean ${{\bf 0}}$ and a separable covariance $\Sigma
= A_0 \otimes I_n$.

To estimate the precision matrix $\Theta=A^{-1}$,  consider the
following regressions, where we regress one variable against all
others: for $X_0 = Z_1 A^{1/2}$ defined in \eqref{eq::addmodel}, the $j$th column of $X_0$ satisfies
\begin{eqnarray}
    \label{eq::regr}
 X_{0,j} & = & X_{0, -j}\beta^{j*} + V_{0, j} = \sum_{k \not=j} X_{0, k} \beta_k^{j*} + V_{0, j},  \; \text{where}\;\; \beta_k^{j*} = -{\theta_{jk}}/{\theta_{jj}}, \\
\label{eqLreggg2} 
    \text{and}     &&            V_{0,j} \sim {\N}(0_n, \sigma_{V_j}^2 I_n)\  \mbox{is independent of}\ X_{0, -j}.
\end{eqnarray}
Here $X_{0,j}$ denotes the $j^{th}$ column of $X_0$, $X_{0,-j}$
denotes the matrix of $X_0$ with its $j^{th}$ column removed.
There exist explicit relations between the regression coefficients, 
error variances and the concentration matrix $\Theta =A^{-1} :=
(\theta_{ij}) \succ 0$: for all $j$, we have $\sigma^2_{V_j} =
1/{\theta_{jj}} >0$, and $\theta_{jk}=  \beta^{j*}_k=  \beta^{k*}_j= 0$  if $X_{0,j}$ is 
independent of $X_{0,k}$ given the other variables.

\section{Models and methods}
\label{sec::envregr}
In this paper, we will adapt the multiple Gaussian regression framework to the subgaussian random 
ensembles.
Note that the relationship \eqref{eq::regr} still holds for the 
general case that $X_0$ consists of dependent sub-gaussian random variables,
but the conditional independence property as in \eqref{eqLreggg2}
only holds when $X_0$ is a Gaussian random ensemble.
In particular, for the sub-gaussian model as considered in the present work, 
instead of independence, we have by the Projection Theorem, for all 
$j$, $\cov(V_{0,j}, X_{0,k}) = 0 \;\forall k \not=j$
cf. Proposition~\ref{prop::projection}.
Denote by $Z$ a  random ensemble where
$Z_{ij}$ are independent 
subgaussian random variables such that 
\ben
\label{eq::Zdef}
\forall i, j \quad
\expct{Z_{ij}} = 0, \; \var(Z_{ij}) = 1, \; \text{ and } \; \; \norm{Z_{ij}}_{\psi_2} \leq K,
\een
for some absolute constant $K$, where for a random variable $X$,
$\norm{X}_{\psi_2} := \mathrm{inf}\{t>0 : \mathbb{E} [\exp(X^2/t^2)]
\le 2\}$.
Let $Z_1, Z_2$ be independent copies of $Z$.
Proposition~\ref{prop::projection},
a reduction from Proposition
2.3. in~\cite{Zhou24}, illuminates this zero correlation condition for
the general sub-gaussian matrix variate model \eqref{eq::addmodel}, as well as the explicit
relations between the regression coefficients $\beta_k^{j*}$, error
variances for $\{V_{0,j}, j \in [n]\}$, and the inverse covariance
$\Theta= (\theta_{ij}) \succ 0$.
\begin{proposition}{\textnormal{[Zhou 2024] (Matrix subgaussian model)}}
\label{prop::projection} 
Let $\{X_{0,j}, j =1, \ldots, m\}$ be column vectors of the data matrix
$X_0 = Z_1 A^{1/2}$ as in~\eqref{eq::addmodel},
where we assume that $A \succ 0$.
Let $\Theta = A^{-1} = (\theta_{ij}) \succ 0$, where $0<
1/{\theta_{jj}} < a_{jj}, \forall j$.
Consider many regressions, where we regress one column vector $X_{0,j}$
against all other column vectors.
Then for each $j \in [m]$,
\ben
\label{eq::subregress}
&& X_{0,j} = \sum_{k \not=j} X_{0,k} \beta_k^{j*} + V_{0,j} \;
\; \text{ where } \quad \beta_k^{j*} = - {\theta_{jk}}/{\theta_{jj}}, \\
\label{eq::subresidual}
&& 
\cov(V_{0,j}, V_{0,j}) = I_n/{\theta_{jj}}  \; \text{ and} \; \cov(V_{0,j}, X_{0,k}) = 0 \; \; \forall k \not=j.
\een
\end{proposition}
\begin{proofof}{Proposition \ref{prop::projection}}
  Let $Y = X_0^T$. Then  $Y = A^{1/2} Z_1^T I_n$ as in~\eqref{eq::addmodel}.
  Then explicit relations between row vectors of $Y$ are characterized in
  Proposition 2.3. where we set $B_0 = A$ and $A_0 = I$, and $\Theta_0
  = A^{-1}$.  Then the results of Proposition 2.3 immediately translate to the relationships among column vectors $\{X_0^j, j =1, \ldots, m\}$.
  This completes the proof.
\end{proofof}

We propose the following procedures for estimating the precision
matrix $\Theta := A^{-1}$ in view of Proposition~\ref{prop::projection}.
Following~\cite{Zhou24}, our method on inverse covariance estimation 
corresponds to the proposal in~\cite{Yuan10} and~\cite{LW12}, only now dropping 
the i.i.d. or Gaussian assumptions, as we use the 
relations~\eqref{eq::subregress} and~\eqref{eq::subresidual}, 
which are valid for the much more general subgaussian matrix-variate 
model~\eqref{eq::addmodel}.
Let $X_{\minus i}$ denote columns of $X$ without $i$ and 
$\hat\tau_B$ be as in \eqref{eq::trBest}.
We obtain $m$ vectors of
$\hat{\beta}^i, i=1, \ldots, m$,  
by solving~\eqref{eq::origin} with the following input:
\ben
\label{eq::gamma}
\hat\Gamma^{(i)} = \onen X_{\minus i}^T  X_{\minus i} - \hat\tau_B
I_{m-1}, \quad \text{and } \; \; \hat\gamma^{(i)} \; = \; \onen  X_{\minus i}^T X_i.
\een
For a chosen penalization parameter $\lambda \geq 0$ and a fixed $b_1
> 0$, consider the following corrected Lasso estimator for the nodewise regression: for each $i \in [m]$,
\begin{eqnarray}
\label{eq::origin}
\label{eq::hatTheta}
\hat \beta^j =
 \argmin_{\beta \in \R^{m-1},  \norm{\beta}_1 \le b_1}
    \frac{1}{2} \beta^T \hat\Gamma^{(j)} \beta - \ip{\hat\gamma^{(j)}, \beta} + \lambda^{j} \|\beta\|_1,
\end{eqnarray}
where $b_1$ is chosen so that $\norm{\beta^{j*}}_1 = \sum_{k=1}^{m-1}
\abs{\beta^{j*}_k} \le b_1$ for all $j \in [m]$.
Clearly, the choice of $b_1$ in~\eqref{eq::hatTheta} will depend on the model class for $\Theta$,
which will be chosen to provide an upper bound on the matrix $\ell_1$
norm $\norm{\Theta}_1$.
In particular, we assume that $\theta_{jj} \le M$ for some absolute
constant $M$ so that the residual error variance $\sigma_{V_j}^2 :=
1/\theta_{jj} \ge 1/M$ is lower bounded for all $j \in [m]$;
cf. Propositions~\ref{prop::projection} and~\ref{prop::XVcorr}. In summary, we have Algorithm 1.

\noindent
\textsc{Algorithm 1: Input ($\hat\Gamma^{(j)}, \hat\gamma^{(j)}$) as 
  in~\eqref{eq::gamma}} \\
\noindent{\textbf{(1)}}
Perform $m$ regressions using \eqref{eq::hatTheta} to obtain vectors of
$\hat{\beta}^{j} \in \R^{m-1}$, $j=1, \ldots, m$, where $\hat{\beta}^j =\{\hat{\beta}_k^j; k \in  
[m], k \not=j\}$, with penalization parameters $\lambda^{j}, b_1 > 0$ to be specified.
More precisely, we can set
\bens
b_1 \asymp \max_{j} (\theta_{jj}
\sum_{i\not=j}^m \abs{\theta_{ij}}) \le M \norm{\Theta}_1,\;
\text{ where} \; 
\norm{\Theta}_{1} = \max_{j}\sum_{i=1}^m \abs{\theta_{ij}};
\eens
\noindent{\textbf{(2)}}
Obtain an estimate for each row of $\Theta$ as follows:  
\bens
\forall j \in [m],
\tilde{\Theta}_{j j} = (\hat\Gamma_{jj} - \hat\Gamma_{j, 
   \minus j} \hat{\beta}^{j})^{-1}, \tilde{\Theta}_{j, \minus j} = -\tilde{\Theta}_{j j} \hat{\beta}^j,
 \eens
 where $\tilde{\Theta}_{j, \minus j}$ denotes the $j^{th}$ row of
 $\Theta$ with diagonal entry $\tilde{\Theta}$ removed; \\
\noindent{\textbf{(3)}}
Project $\tilde{\Theta}$ onto the space $\SM$ of $m \times m$ symmetric matrices:
\[
\hat\Theta = \argmin_{\Theta \in \SM} \shonenorm{\Theta - \tilde{\Theta}}.
\]
\subsection{The errors-in-variables regression}
\label{sec::regress}
Let $\{X_{0,j}, j =1, \ldots, m\}$ be column vectors of the data matrix $X_0 =
Z_1 A_0^{1/2}$ as in~\eqref{eq::addmodel}.
Now the $j$th column of $X$ in \eqref{eq::addmodel}  can be written
as,
\ben
\label{eq::additive}
&& \forall j \in [m], \; X_j  =  X_{0,-j} \beta^{j*} + V_{0,j} + W_j =: X_{0,-j} \beta^{j*}  + \ve_j,   \\
\nonumber
&& \text{where we observe  $X_j$ and } \; X_{-j} = X_{0, -j} + W_{-j},
\; \text{instead of } \quad X_{0,-j}
\een
Moreover, the noise vector $\ve_j :=   V_{0,j} + W_j$ in \eqref{eq::additive} may no longer be independent of $\{X_{0, 
  i}; i \neq j\}\ (i \in [m])$, but they remain uncorrelated; cf \eqref{eq::subresidual}.
Moreover, the components of $\ve_j$ are correlated due to $W_j$.
However, this model still fits in the errors-in-variables
framework despite the complication due to the dependencies between
components of $\ve$ as we will show in Section~\ref{sec:graphical}.
A major part of our theory relies on untangling 
the dependencies among $\ve_j$ and $X_{0, -j}$. 
Proposition~\ref{prop::XVcorr}
describes the correlation between the
residual error vector $ V_{0, j}$  for each $j \in [m]$ and the $m-1$ column
vectors $X_{0, i}, i \in [m], i\not=j$ in the regression
function~\eqref{eq::subregress}.
This result is new to the best of our knowledge. See~\cite{Zhou24} for 
results related to matrix-variate data with missing values. 
\begin{proposition}
\label{prop::XVcorr}
  Suppose $n =\Omega(\log m)$.
  Let $K := \max_{i,j} \norm{Z_{ij}}_{\psi_2}$.
  Let $X_0 = Z A^{1/2}$,  where $A^{1/2}$ denotes the unique square 
  root of the positive definite matrix $A$ and $Z$ consists of independent entries satisfying~\eqref{eq::Zdef}.
  Denote by $A^{1/2} = [d_1, \ldots, d_m]$, where 
$d_1, \ldots, d_m \in \R^m$.
  Denote by $X^{\ell}_0 = (X^{\ell}_{0,1}, X^{\ell}_{0,2}, \ldots, 
X^{\ell}_{0,m}), \ell \in [n]$ any row vector from the data matrix $X_0$.
 Let $V_{0, j} = (V^{1}_{0, j}, \ldots, V^{n}_{0, j})$ be as defined
 in~\eqref{eq::subregress}.
 Let $C, C', C_0$ be absolute constants.
Then for all $\ell \in [n], j \in [m]$,
\bens 
\norm{X^{\ell}_{0, j}}_{\psi_2}
& \le & C K \twonorm{d_j} =  C  a^{1/2}_{jj} \;
\text{and} \; \norm{V^{\ell}_{0, j}}_{\psi_2} \le
C' K \sigma_{V_j}, \; \;\text{ where} \; \; 
\sigma_{V_j}^2 = 1/{\theta_{jj}}.
\eens
With probability at least $1-\inv{m^2}$, $\max_{i \not= j} \inv{n} \ip{X_{0, i}, V_{0, j}} \le C_0
  K^2 \sigma_{V_j} \sqrt{a_{\max}} \sqrt{\log m/n}.$
\end{proposition}

\subsection{Notation and organization}
For a matrix $A = (a_{ij})_{1\le i,j\le m}$, we use $\twonorm{A}$ to denote its operator norm and
let $\norm{A}_{\max} =\max_{i,j} |a_{ij}|$ denote  the entry-wise max norm.
Let $\norm{A}_{1} = \max_{j}\sum_{i=1}^m\abs{a_{ij}}$ denote the matrix $\ell_1$ norm.
For a square matrix $A$, let ${\rm tr}(A)$ be the trace of $A$, $\diag(A)$ be a diagonal matrix with the same diagonal as
$A$, and $\kappa(A)$ denote the condition number of $A$.
Let $I_n$ be the $n$ by $n$ identity matrix.
For a number $a$,  $\mathrm{sign}(a)=1$, $-1$, and $0$ when $a>0$, $a<0$, and $a=0$, respectively.
For two numbers $a$ and $b$, $a \wedge b
:= \min(a, b)$ and $a \vee b := \max(a, b)$.
We use $a =O(b)$ or $b=\Omega(a)$ if $a \le Cb$ for some positive absolute
constant $C$ which is independent of $n, m$ or sparsity parameters.
We write $a \asymp b$ if $ca \le b \le Ca$ for some positive absolute
constants $c,C$. The absolute constants $C, C_1, c, c_1, \ldots$ may change line by
line.

\noindent{\bf  Organization}
The rest of this paper is organized as follows. 
Section~\ref{sec:theory} presents theoretical results on inverse 
covariance estimation, where we develop convergence bounds in 
Theorem~\ref{coro::Theta} for the matrix-variate subgaussian
model~\eqref{eq::addmodel}.
We prove Proposition \ref{prop::XVcorr} and Lemma~\ref{lemma::low-noise} in Section~\ref{sec::proofofmain}.
Proof of Theorem~\ref{coro::Theta} appears in Section \ref{sec:thm1}.
In Section~\ref{sec::conclude}, we conclude.

\section{Theoretical properties}
\label{sec:theory}
In this section, we establish statistical convergence and the model
selection consistency for the inverse covariance estimation problem
for the matrix-variate subgaussian distribution~\eqref{eq::addmodel}
with Kronecker sum covariance model.  We require the following assumptions.\\
\noindent{\bf (A2)}
The minimal eigenvalue $\lambda_{\min}(A)$
of the covariance matrix $A$ is bounded: $\lambda_{\min}(A) > 0$, and 
the condition number $\kappa(A) =O(\sqrt{{n}/{\log m}})$.

\noindent{\bf (A3)}
Suppose $\tau_B =
O(\lambda_{\max}(A))$. The covariance matrix $B$ satisfies ${\fnorm{B}^2}/{\twonorm{B}^2}
=\Omega(\log m)$ and
\bens
{\tr(B)}/{\twonorm{B}}  = \Omega(d \log (m/d)) \; \text{ where}  \; \;
d \asymp n/\log m
\eens
Note that Assumption (A2) requires that $A \succ 0$ and the condition
number $\kappa(A)$ is upper bounded.
(A3) requires that $\tr(B)$ is both upper and lower bounded.
The smallest possible value for the stable rank ${\fnorm{B}^2}/{\twonorm{B}^2}$ and 
the effective rank ${\tr(B)}/{\twonorm{B}}$ is $1$ when $B$ is a 
rank-one matrix, while these quantities tend to have larger values
when the maximum eigenvalue of $B$ becomes closer to the remaining
eigenvalues.
See Definition~\ref{def::lowRE}, Corollary 25, Theorem 26 and Lemma 
15 in~\cite{RZ15}, where conditions (A1) -(A3) are stated for the Lower 
and Upper-RE conditions to hold for the modified Gram matrix 
$\hat\Gamma$~\eqref{eq::hatGamma} (and hence $\hat\Gamma^{(j)}$ as in~\eqref{eq::gamma}). 
Theorem \ref{coro::Theta} shows the statistical consistency of $\hat\Theta$ in the operator norm. 

\begin{theorem}
\label{coro::Theta}
Consider~\eqref{eq::addmodel}.
Suppose conditions (A1)-(A3) hold.
Suppose the columns of $\Theta$ is $d$-sparse, i.e., the number of
nonzero entries on each column in $\Theta$ is upper bounded by $d$,
which satisfies $d=O(n/ \log m)$. Suppose the condition number
$\kappa(\Theta) <\infty$.
Let $\hat{\beta}^i$ be an optimal solution to the nodewise regression
\eqref{eq::origin} with $b_1 := b_0 \sqrt{d}$, where $b_0$ satisfies
$\phi b_0^2 \le \twonorm{\beta^{i*}}^2 \le b_0^2$ for some $0< \phi
<1$.
Let $D'_0 = {\twonorm{B}}^{1/2} + a_{\max}^{1/2}$ and $C_0$ be some
positive absolute constant. Let
\ben
  \label{eq::lambda}
  \forall i, \lambda^{i}  & \ge & 4  \psi_i\sqrt{{\log m}/{n}} \;  \text{ where }  \;
\psi_i := C_0 D_0' K^2  \left(\tau_B^{+/2} \twonorm{\beta^{i*}} + \sigma_{V_i}  \right),
\een
where $\sigma_{V_i}^2 := A_{ii}-A_{i,-i} A_{-i,-i}^{-1} A_{-i,i}$ and $K$ is as in~\eqref{eq::Zdef}.
Then,
\begin{equation}
\label{eq:bound:thm1}
\twonorm{\hat\Theta - \Theta}
=O_P\left(d (\max_i \lambda^{i} ) K^2(A)/{\lambda_{\min}^2(A)}  
\right). 
\end{equation}
Moreover, suppose that $\tau_B \le \twonorm{B} = O(\lambda_{\max}(A))$ and
$ \lambda^{i}  \asymp \psi_i\sqrt{\frac{\log m}{n}}$ for each $i$. Then
\begin{equation}
\label{eq:rel:bd}
{\twonorm{\hat{\Theta}-\Theta}}/{\|\Theta\|_2} = O_P\left(d \kappa^4(A)
  \sqrt{{\log m}/{n}}\right).
\end{equation}
\end{theorem}
Using similar arguments,  we can show consistency results of $\Phi = B^{-1}$.   See Theorem \ref{coro::Omega} of the Supplementary materials for details.

\begin{remark}
\label{rem::factor}
In \eqref{eq::lambda},  since $\beta^{i*}= -(A_{ii}-A_{i,-i} A_{-i,-i}^{-1} A_{-i,i}) \Theta_{-i,i}$,
\bens
 \|\beta^{i*}\|_2 \le a_{\max}  \|\Theta_{-i,i}\|_2 \le a_{\max}
 \lambda_{\max} (\Theta) =
\frac{a_{\max}}{\lambda_{\min}(A)},
\eens
where $a_{\max} := \max_{i} A_{ii} \ge 1$ by (A1). Together with
$\sigma_{V_i}^2 < a_{ii} \le a_{\max}$,
cf. Proposition~\ref{prop::projection} and  $\tau_B^{+/2} =
O(\sqrt{\tau_B})$,
we have $\lambda^{i} \asymp \twonorm{B} \kappa(A) \sqrt{{\log m}/{n}}$, and
hence for $\lambda_{\min}(A) = \twonorm{\Theta}^{-1}$
\bens
\twonorm{\hat\Theta - \Theta}/\twonorm{\Theta} = 
O_p( \frac{d \twonorm{B} \kappa^3(A)}{\lambda_{\min}(A)}   \sqrt{{\log m}/{n}}).
\eens
Thus~\eqref{eq:rel:bd} holds so long as $\twonorm{B} =
O(\lambda_{\max}(A))$. When the noise in  $W$ is negligible in the sense that $
\tau_B$ (and hence $\|B\|_2$ in view of (A3)) is close to zero, then
we have $\lambda^{i} \asymp a_{\max} \sqrt{{\log m}/{n}}$, and
\eqref{eq:bound:thm1} reduces to
$\twonorm{\hat\Theta - \Theta}/\twonorm{\Theta}  =
O_P\left(d \kappa^2(A) \frac{a_{\max}}{\lambda_{\min}(A)} \sqrt{{\log m}/{n}}\right)$,
 recovering the rate of convergence obtained by the Lasso when $X_0$ is free
of measurement error.
\end{remark}

The bound \eqref{eq:bound:thm1} is analogous to that of Corollary 5 in
\citet{LW12}, where the measurement error $W$ has independent row
vectors with a known $\Sigma_W := \onen \E W^T W$.
In the present work, assuming (A1) allows us to estimate $\E[W^T W]$, given the 
knowledge that $W$ is composed of independent column vectors.
In contrast, existing work needs to assume that the covariance matrix 
$\Sigma_W := \onen \E [W^T W]$ of the independent row vectors of $W$
or its functionals are either known a priori, or can be estimated from
a dataset independent of $X$, or from replicated $X$ measuring the
same $X_0$; see for example \cite{carr:rupp:2006}, \cite{RT10,RT13},
\cite{LW12} and \cite{BRT14}.
As described in supplementary Section \ref{sec:replicate},  we can use
replicated measurements of $X_0$ (replicates) to obtain an estimator $\hat{B}$ for covariance $B$. Then,
we can set $\hat{\tau}_B = \tr(\hat{B})/n$ and hence we do not need
(A1).
Although the model we consider is different from those in the
literature, the identifiability issue, which arises from the fact that
we observe the data under an additive error model, is common. Such
repeated measurements are not always available or costly to obtain in
practice~\citep{carr:rupp:2006}.

\subsection{Proof sketch of Theorem \ref{coro::Theta}}
\label{sec:graphical}
Lemmas~\ref{lemma::low-noise} and~\ref{lemma::AD} are essential to prove 
Theorem~\ref{coro::Theta}. 
It is analogous to  Theorem 12 and Corollary 13 of \citet{RZ15}.
However, the key difference between the well-established errors-in-variables analysis
in~\cite{RZ15} and that in the present work is the noise vector $\ve_i = (\ve_{i,1}, \ldots \ve_{i,n})$ that appears in the linear model~\eqref{eq::additive},
\ben
\nonumber
\ve_i & := &  V_{0,i} + W_i, \; \; \text{ where } \;\;
 V_{0,i}  \; \mbox{is isotropic and } \; \cov(V_{0, i}, X_{0, j}) = 0 \; \forall \; j
 \neq i, \\
 \label{eq::WXindep}
 \text{ and } &&  W_i \; \;  \mbox{ is independent of} \; \; X_{0} \; \; \text{ and 
} \; V_{0} \; \; \; (i=1,\ldots,m).
\een
Compared to the Gaussian regression model, it remains true that
$(W_i)_{i \in [m]}$ is independent of  $X_{0}$ and correspondingly
$V_0$ ; However, the noise components $(\ve_{i,1}, \ldots \ve_{i,n}), (i \in  [m])$ are now correlated due to the correlation in $W_i$ and
they may no longer be independent of $\{X_{0,  j}; j \neq i\}$.
Let
\ben
\label{eq::defineD0}
& & D_0 = \sqrt{\tau_B} + \sqrt{a_{\max}}, \; \; \;
 D_0' =
\twonorm{B}^{1/2} + \sqrt{a_{\max}}, \; \; \tau_B^+ :=
(\tau_B^{+/2})^2, \\
\label{eq::defineDtau}
&&
\text{where} \; \;  \tau_B^{+/2}
:=   \sqrt{\tau_B} + {D_{\ora}}/{\sqrt{m}}, \; \;
D_{\ora} \;  = \; 2 (\twonorm{A}^{1/2} + \twonorm{B}^{1/2}).
\een
\begin{lemma}
\label{lemma::low-noise}
Suppose (A1) holds.
Suppose that $m \ge 16$ and
${\fnorm{B}^2}/{\twonorm{B}^2} \ge \log m.$
Let $\hat\tau_B$ be defined as in~\eqref{eq::trBest}.
Let $\hat\Gamma^{(i)}$ and $\hat\gamma^{(i)}$ be as defined in Algorithm 1.
On event $\B_0$, where $\prob{\B_0} \ge 1- 16/m^2$, we have for $\hat \beta^i$ as in~\eqref{eq::hatTheta} and 
some absolute constant $C_0$, 
\ben
\label{eq::psioracle}
\forall i, \quad
\norm{\hat\gamma^{(i)}
 - \hat\Gamma^{(i)} \beta^{i*}}_{\infty}
& \le &
\psi_i \sqrt{\frac{\log m}{n}}, \text{ where} \\
\nonumber
\psi_i := C_0 K^2 D_0'\left(\tau_B^{+/2} \twonorm{\beta^{i*}} + \sigma_{V_i}
\right) & \le & C_0 K^2 D_0' (\tau_B^{+/2} \kappa(A) + \sigma_{V} ).
\een
\end{lemma}
\begin{lemma}
  \label{lemma::AD}
Suppose the conditions in Lemma~\ref{lemma::low-noise} hold.
Let $\Delta := \hat\Gamma -A := \onen X^TX - \hat\tau_B I_{m} -A.$
Let $D_1 :=   \frac{\fnorm{A}}{\sqrt{m}} +
\frac{\fnorm{B}}{\sqrt{n}}$ and $\psi_0 = C_0 K^2 (\twonorm{B}+  a_{\max})$.
Then with probability at least $1- 12/m^2$,
\ben
\label{eq::defineAD}
\maxnorm{\Delta}
 \le  8 \tau_B^{1/2} (D'_0 + a_{\max})  r_{m,n} + 4 C_0 D_1 r_{m,m}
  \le    12 C K^2 \vp \sqrt{{\log  m}/{n}}  \asymp \psi_0\sqrt{{\log m}/{n}}
  \een
\end{lemma}
Such error bounds are important in the regular Lasso error bounds for 
nodewise regression~\citep{ZRXB11,ZhouTH23} and the corrected Lasso for 
errors-in-variables regression and Gaussian graphical 
models~\citep{LW12,RZ15}.
Our theoretical results show that $\hat{\Theta}$ consistently estimates $\Theta$ in the
spectral norm under suitable conditions.
However, $\hat{\Theta}$  is not necessarily positive-semidefinite.
One can obtain a consistent and positive-semidefinite estimator
by considering an additional estimation procedure, which is described
in Section \ref{sec:add:proc} of the Supplementary materials.

\subsection{Proof of  Theorem~\ref{coro::Theta}}
\label{sec:thm1}

We consider the following variation of the baseline $\RE$ condition.
\begin{definition}{\textnormal{(Lower-$\RE$ condition)~\citep{LW12}}}
\label{def::lowRE}
The $m \times m$ matrix $\Gamma$ satisfies a Lower-$\RE$ condition with curvature
$\alpha >0$ and tolerance $\tau > 0$ if
\bens
\theta^T \Gamma \theta \ge
\alpha \twonorm{\theta}^2 - \tau \|\theta\|_1^2, \; \;  \forall \theta \in \R^m.
\eens
\end{definition}
As $\alpha$ becomes smaller, or as $\tau$ becomes larger, the
Lower-$\RE$ condition is easier to be satisfied.
\begin{proofof}{Theorem~\ref{coro::Theta}}
  Lemma~\ref{lemma::low-noise} and the supplementary Lemma~\ref{lemma::AD}
  show that the following conditions hold with probability at least  $1- 20 /m^2 - 6/m^3$:
  \ben
  \nonumber
\norm{\hat\gamma^{(i)}- \hat\Gamma^{(i)} \beta^{i*}}_{\infty} & \le &  C_\psi \sqrt{{\log
    m}/{n}} < \lambda^{i}/4, \\
\label{eq::gammalocal}
\maxnorm{\hat\Gamma -A } & := & \maxnorm{\onen X^TX - \hat\tau_B I_{m}
  -A} = O\left(C_\psi \sqrt{{\log m}/{n}}\right),
\een
where $C_\psi := K^2 C_0 D_0'  \left(\tau_B^{+/2} \kappa(A)+ \sqrt{a_{\max}} \right)$.
~\cite{RZ15} show that the lower-$\RE$ condition holds with high
probability with curvature $\alpha = 5\lambda_{\min}(A)/8$ and
tolerance $\tau =\frac{\lambda_{\min}(A)}{2s_0} \asymp
\frac{\alpha}{s_0}$ uniformly over all $\hat\Gamma^{(i)}$.
The theorem then follows from Corollary 5 in~\cite{LW12}, so long as 
we can show that condition~\eqref{eq::taumain} holds for $\lambda^{(i)} \ge 4 \psi_i\sqrt{\frac{\log m}{n}}$, 
where  the parameter $\psi_i$ is as defined in~\eqref{eq::lambda}:
\ben
\label{eq::taumain}
\sqrt{d} \tau \le \min \left\{\frac{\alpha}{32 \sqrt{d}}, \frac{\lambda^{(i)}}{4b_0} \right\}.
\een
Condition \eqref{eq::taumain} can be easily checked using the 
definition of $s_0 =O(n/\log m)$ as in \eqref{eq::s0cond} and the lower bound of 
$\lambda^{(i)}$ as already shown in the proof of Theorem 3 of~\cite{RZ15}. 
Hence, the theorem indeed follows from Corollary 5 in~\cite{LW12}. 
\end{proofof}

\subsection{Related work and discussions}
\label{sec::related}
Variants of the linear errors-in-variables models in the high
dimensional setting has been considered in recent
work~\citep{RT10,LW12,RT13,CC13,SFT14,BRT14}, where oblivion in the
covariance structure for row or columns of $W$, and a general
dependency condition in the single data matrix $X$ are not
simultaneously allowed.
An estimator similar to~\eqref{eq::origin} was considered by~\cite{LW12}, 
which is a variation of the  Lasso \citep{Tib96} or the Basis 
Pursuit~\citep{Chen:Dono:Saun:1998} estimator.
The non-convex optimization function in \eqref{eq::origin} can be
solved efficiently using the composite gradient descent algorithm as
described in supplementary Section \ref{compsec}; see \citet{ANW12}, \citet{LW12}, and
\citet{RZ15} for details.
For the related Conic programming estimators (and related Dantzig-type
methods) for estimating EIV regression coefficients,
see~\cite{RT10,RT13}, \cite{BRT14}, and \cite{RZ15}.

For matrix-variate data with two-way dependencies, 
 prior work depended on a large number of replicated data to obtain certain convergence guarantees,
 even when the data is observed in full and free of measurement error;
 see for example~\cite{Dut99}, \cite{WJS08}, and~\cite{THZ13}.
A recent line of work on high dimensional matrix variate models \citep{KLLZ13,Zhou14a,GZH19,ZG24,Zhou24}
have focused on the design of estimators and computationally efficient algorithms
while establishing sharp rates of convergence on the inverse
covariance estimation.
~\cite{Zhou24} developed multiple regression methods for estimating 
the inverse covariance matrices in a subgaussian matrix variate model, 
where the matrix variate data is observed through a random mask U.
See also~\cite{Efr09}, \cite{AT10}, and \cite{HSZ15} for related 
models and applications of the Kronecker product (inverse) covariance 
model.

Among these models, the Kronecker sum provides a 
 covariance or precision matrix which is sparser than the  Kronecker product. 
In particular, the models considered in~\cite{KLLZ13, GZH19,ZG24}
use the Kronecker Sum to model the inverse covariance matrix directly.
In contrast, in the current paper, motivated by the
errors-in-variables literature, we consider data matrix $X$ generated
from~\eqref{eq::addmodel} with the Kronecker sum as the covariance
matrix for $X$ as in~\eqref{eq:model_intro0},
while simultaneously exploiting the multiple regression
framework to estimate the inverse of each component in the KS-based
covariance matrix; cf.~\eqref{eq::additive}.
Previously,~\cite{KSSA:2017} considered a slightly more
sophisticated model based on the sum of Kronecker Products to model
time-varying component $W = [w^1, \ldots, w^m]$ that consists of
column vectors whose covariance changes smoothly over time; that is,
we allow $E (w^t \otimes w^t) = B(t)$ where $B(t)$ is smoothly varying
over time.
The methodology in~\cite{KSSA:2017} essentially follows
from~\cite{ZLW08a}, using the graphical Lasso with the regularized
kernel estimator of the sample covariance at time $t$ to estimate
$B(t)^{-1}$ for all $t$. In the current paper, we use the
multiple regression based methods, but with fixed $B(t) \forall t$.
We refer to~\cite{GZH19,ZG24} for other applications of the Kronecker 
Sum models and further references.

\section{Proof of Proposition \ref{prop::XVcorr} and Lemmas
  \ref{lemma::low-noise} and~\ref{lemma::AD} }
\label{sec::proofofmain}
The sub-exponential (or $\psi_1$) 
norm of random variable $S$, denoted by $\norm{S}_{\psi_1}$, is defined as 
$\norm{S}_{\psi_1} = \inf\{t > 0\; : \; \E \exp(\abs{S}/t) \le 2 \}.$
First, we state the Bernstein’s inequality.
\begin{theorem}{\textnormal{(Bernstein’s Inequality)~\cite{Vers18}}}
  \label{thm::berns}
 Let $X_1, \ldots, X_{n}$ be independent, mean zero, sub-exponential
 random variables. Let $L =  \max_{i} \norm{X_i}_{\psi_1}$.
 Then, for every $t \ge 0$, we have
 \bens 
 \prob{\abs{\inv{n} \sum_{i=1}^n X_i} >t} \leq 
 2 \exp\big(-c \min(\frac{t^2}{L^2}, \frac{t}{L} ) n\big).
\eens
\end{theorem}

\begin{lemma}{\textnormal{(Proposition 2.7.1) ~\cite{Vers18}}}
\label{lemma::product}(Product of sub-gaussians is sub-exponential).
Let $X$ and $V$ be sub-gaussian random variables. Then $XV$ is sub-exponential. Moreover,
\bens
\norm{XV}_{\psi_1} \le  \norm{X}_{\psi_2} \norm{V}_{\psi_2} 
\eens
\end{lemma}

\begin{proofof}{Proposition \ref{prop::XVcorr}}
  Denote by $\Theta_{j \cdot}$ the $j^{th}$ row vector of $\Theta$.
Then by definition, $X^{\ell T}_0 = (Z^{\ell}_{1},
Z^{\ell}_{2},\ldots, Z^{\ell}_{m})^T A^{1/2}$, where $Z^{\ell}$
denotes the $\ell^{th}$ row vector of $Z$. Hence
$X^{\ell}_{0, k} =  e_k^T X^{\ell}_0 = e_k^T  A^{1/2}
(Z^{\ell}_{1}, Z^{\ell}_{2}, \ldots, Z^{\ell}_{m})$.
Hence for all $\ell \in [n], k \in [m]$,
\bens
\norm{X^{\ell}_{0, k}}^2_{\psi_2} & = &
\norm{d_k^T (Z^{\ell}_{1}, Z^{\ell}_{2}, \ldots, Z^{\ell}_{m})
}^2_{\psi_2} \le  C K^2 \twonorm{d_k}^2 =  C K^2 a_{kk}
\eens
by Proposition 2.6.1~\cite{Vers18}.
On the other hand, by the linear relationship, we have for all $j = 1, \ldots, 
m$ and for any $\ell \in [n]$,
\bens
V^{\ell}_{0,j}
& = & X^{\ell}_{0,j} -\sum_{k \not=j} X^{\ell}_{0,k} \beta_k^{j*} \\
& = &
X^{\ell}_{0,j} + \sum_{k \not=j} X^{\ell}_{0,k}
\frac{\theta_{jk}}{\theta_{jj}}
=  \inv{\theta_{jj}} \Theta_{j \cdot}^T X^{\ell}_{0}  =  \inv{\theta_{jj}} \Theta_{j \cdot}^T  A^{1/2} Z^{\ell}
\eens
Then for all $\ell \in [n], j \in [m]$, and $\sigma_{V_j}^2 = 1/{\theta_{jj}}$,
\bens
\norm{V^{\ell}_{0, j}}^2_{\psi_2}
& \le &
\frac{C} {\theta^2_{jj}} \twonorm{ A^{1/2} \Theta_{j \cdot}}^2
\max_{j, \ell}\norm{Z_{j}^{\ell}}^2_{\psi_2}   \le  
\frac{C K^2}{\theta^2_{jj}} \Theta_{j \cdot}^T A \Theta_{j \cdot} \\
& = & \frac{C K^2}{\theta^2_{jj}} (\Theta A \Theta)_{jj}= C K^2/{\theta_{jj}} =: C K^2 \sigma_{V_j}^2.
\eens
  It remains to prove the large deviation inequality.
  W.l.o.g., we fix $i =1$ and $j = 2$. Then
 by \eqref{eq::subresidual}, 
\bens 
\forall k \in [n], 
\quad  \E (X_{0, 1}^{k} V_{0, 2}^{k}) & =  & \E (X_{0, 1}^k) \E(V_{0,
  2}^k) = 0, \text{ and }
\onen \ip{X_{0, 1}, V_{0, 2}} = \onen \sum_{\ell=1}^n X^{\ell}_{0, 1} V^{\ell}_{0, 2},
 \eens
  where $X_{0, 1} =(X^{1}_{0, 1}, \ldots, X^{n}_{0, 1})$  (resp.
$V_{0, 2} =(V^{1}_{0, 2}, \ldots,  V^{n}_{0, 2})$)  consists of
independent, mean-zero, subgaussian components whose $\psi_2$ norms are uniformly
bounded as above.
Now by Lemma~\ref{lemma::product}, we have for some absolute
constant $C$, for all $\ell \in [n]$,
\bens
\norm{X^{\ell}_{0, i} V^{\ell}_{0, j}}_{\psi_1} \le 
\norm{X^{\ell}_{0, i}}_{\psi_2} \norm{V^{\ell}_{0, j}}_{\psi_2} s
\le C K^2 \sqrt{a_{ii}} \sigma_{V_j} \le C K^2 \sqrt{a_{\max} } \sigma_{V_{j}}.
\eens
Then for $t_{ij} = C_1 K^2  \sqrt{a_{ii}} \sigma_{V_j} \sqrt{\log m/n}
\asymp \sqrt{\log m/n} L_j$, where $L_j = \max_{i} \max_{\ell}
\norm{X^{\ell}_{0, i} V^{\ell}_{0, j}}_{\psi_1} \le  C'
\sqrt{a_{\max}} K^2 \sigma_{V_j}$, by Bernstein’s Inequality,
cf. Theorem~\ref{thm::berns} (Corollary 2.8.3~\cite{Vers18}),
  \bens 
\prob{\exists i \not= j, \abs{\onen \sum_{\ell=1}^n X^{\ell}_{0, i} V^{\ell}_{0, j} }> t_{ij} }
& \le &  2 n(n-1) \exp\big(-c \min(\frac{t_{ij}^2}{L_j^2}, \frac{t_{ij}}{L_j}) n\big) \\
& =&  2 \exp\big(-c' \min(\log m, \sqrt{n \log   m}) \big)
\eens
This completes the proof.
\end{proofof}

Denote by $\bar\beta^j \in \R^m$ the zero-extended $\beta^{j*}$ in $\R^m$ such that 
$\bar\beta^j_i = \beta^{j*}_i \; \; \forall i \ne j$ and $\bar\beta^j_j = 0$. 
The large deviation bounds in Lemmas 5 and 11 of \citet{RZ15} are the key results in proving Lemmas \ref{lemma::AD} and \ref{lemma::low-noise}. 
Let $C_0$ be an absolute constant appropriately chosen.
Throughout this proof, we denote by:
\bens
r_{m,n} =  C_0 K^2 \sqrt{\frac{\log m}{n}} \; \; \text{ and } \; \
r_{m,m} =  2 C_0 K^2 \sqrt{\frac{\log m}{mn}}.
\eens
Recall $\fnorm{B} /\sqrt{n}\le \sqrt{\tau_B}\twonorm{B}^{1/2}$.
On event $\B_{10}$, which holds with probability at least $1 - 2
/m^2$, cf. Lemma 11~\citep{RZ15}, we have
\ben
\label{eq::B10a}
\onen
\maxnorm{X_0^T W} & \le &  C_0 K^2  \sqrt{\tau_B a_{\max} } 
\sqrt{{\log  m}/{n}}, \text{ and } \\
\label{eq::B10b}
\onen \maxnorm{W^T W - \tr(B) I_m}& \le & 
C_0 K^2  ({\fnorm{B}}/{\sqrt{n}})  r_{m,n} \le \sqrt{\tau_B} \twonorm{B}^{1/2} r_{m,n}
\een
\begin{proofof}{Lemma \ref{lemma::low-noise}}
A careful examination of the proof for Theorem~1 in~\cite{RZ15} shows
that the key new component is in analyzing the following term
$\hat\gamma^{(i)}$ for all $i$.
First notice that for all $i$,
\bens
n \hat\gamma^{(i)}
& = &
{X_{\minus i}^T X_i }  = 
(X_{0, \minus i}^T +   W_{\minus i}^T)(X_{0, \minus i} \beta^{i*} +
V_{0, i} + W_i) \\
& = &
(X_{0, \minus i}^T X_{0, \minus i} \beta^{i*} +  W_{\minus i}^T X_{0, \minus i} \beta^{i*} +
(X_{0, \minus i}^T +   W_{\minus i}^T)(V_{0, i} + W_i) 
\eens
where
$\ve_i :=  V_{0, i} + W_i$, and
\bens
\hat\Gamma^{(i)} \beta^{i*}
& = &{\onen(X_{\minus i}^T X_{\minus i} - \hat\tr(B) I_{m-1}) \beta^{i*}} \\
& = &
\onen (X_{0, \minus i}^T X_{0, \minus i} \beta^{i*}
 + W_{\minus i}^T X_{0,\minus i} \beta^{i*}+
X_{0, \minus i}^T W_{\minus i} \beta^{i*}
+ W_{\minus i}^T W_{\minus i}  \beta^{i*} - \hat\tr(B) I_{m-1} \beta^{i*}) 
\eens
Thus
\bens
&&\norm{\hat\gamma^{(i)} - \hat\Gamma^{(i)} \beta^{i*}}_{\infty}
 \le
\onen\norm{ X_{0, \minus i}^T \ve_i +  W_{\minus i}^T \ve_i
- (X_{0, \minus i}^T W_{\minus i} + W_{\minus i}^T W_{\minus i} - \hat\tr(B) I_{m-1})\beta^{i*}}_{\infty}\\
 & \le &
\onen\norm{ X_{0, \minus i}^T \ve_i +  W_{\minus i}^T \ve_i }_{\infty}
+ \onen\norm{(Z^T B  Z- \tr(B) I_{m}) \bar\beta^i}_{\infty} + \onen \norm{X_0^T W
     \bar\beta^i}_{\infty} \\
&& + \onen \abs{\hat\tr(B) - \tr(B)} \norm{\beta^{i*}}_{\infty} =:  U_1+  U_2 + U_3 + U_4.
\eens
Recall for all $\ell \in [n], j \in [m]$,
$\norm{V^{\ell}_{0, j}}_{\psi_2}
\le C' K \sigma_{V_j}$, \text{  where} \; $\sigma_{V_j}^2 = 1/\theta_{jj}.$
Thus we have by independence of 
$W_0$ and $V_0$, on event $\B_4$,  for every $i \in [n]$,
\bens
\text{ and } \; \;
\max_{j \not= i} \onen \ip{W_j, V_{0, i}} & \le &
 C_0 K^2 \sigma_{V_i} \sqrt{\tau_B}\sqrt{{\log m}/{n}}
= r_{m,n}\sigma_{V_i} \sqrt{\tau_B},
\eens
following the proof of Lemma 11 of \citet{RZ15}. Moreover, we obtain using Proposition~\ref{prop::XVcorr}, on event
$B_8$, which holds with probability at least $1-\inv{n^2}$,
\bens
\max_{j \not=i} \onen \ip{X_{0, j}, V_{0,i}} &  \le &
C_0 K^2 \sigma_{V_i}
\sqrt{ a_{\max}} \sqrt{\frac{\log m}{n}} =
\sigma_{V_i} r_{m,n}\sqrt{ a_{\max}}.
\eens
Here $C_0$ is adjusted so that the error probability hold for $D_0 = \sqrt{\tau_B} + a_{\max}^{1/2}$.
In particular, $U_4$ is bounded using  Lemma 5 \citet{RZ15}, 
which ensures that on event $\B_6$, which holds with probability at least $1-3/m^3$
\bens
\abs{\hat\tr(B) - \tr(B)} /n \le D_1 r_{m,m} \quad \text{ and
  hence}\quad \onen \abs{\hat\tr(B) - \tr(B)}
\norm{\beta^{i*}}_{\infty} \le D_1 r_{m,m} \norm{\beta^{i*}}_{\infty}
\eens
Hence, by  Lemma 11 of \citet{RZ15},   
on event $\B_{10}$, for $W \sim B^{1/2} Z$, which is independent of 
$X_{0}$ and $V_{0}$~\eqref{eq::WXindep},
\bens
\max_{i\not=j} \onen \ip{X_{0, j}, W_i}  & \le &
 \onen\maxnorm{X_0^T W} \le   \sqrt{\tau_B} \sqrt{ a_{\max}} r_{m,n}, \\
 \max_{i \not= j} \onen \ip{W_i, W_j} & \le &\onen\maxnorm{W^T W -\tr(B)I_m} 
\le  \sqrt{\tau_B} \twonorm{B}^{1/2} r_{m,n} 
\eens
On event $\B_5$ for  $D'_0 :=\sqrt{\twonorm{B}} + a_{\max}^{1/2}$, for all $j$,
\bens
\lefteqn{U_2 +  U_3 =
  \onen\norm{(W^T W- \tr(B) I_{m}) \bar\beta^i}_{\infty} + \onen \norm{X_0^T W \bar \beta^j}_{\infty}} \\
& \le & r_{m,n} \twonorm{\beta^{j*}}
\left(\frac{\fnorm{B}}{\sqrt{n}} + \sqrt{\tau_B}
  a^{1/2}_{\max}\right) \le r_{m,n} \twonorm{\beta^{j*}}  \tau_B^{1/2} D_0'.
\eens
Denote by $\B_0 := \B_4 \cap \B_5 \cap \B_6 \cap B_8 \cap \B_{10}$.
Suppose that event $\B_0$ holds.
By Lemmas 5 and 11 of \citet{RZ15}, under (A1) and $D_1$ defined therein, and
Proposition~\ref{prop::XVcorr} in the present paper, on event $\B_0$, for all $i$,
\ben
\nonumber
\norm{\hat\gamma^{(i)} - \hat\Gamma^{(i)} \beta^{i*}}_{\infty}
& \le &
\nonumber
 r_{m,n} \sigma_{V_i} D_0 + D_0' \tau_B^{1/2} r_{m,n}
( \twonorm{\beta^{i*}} + 1)
+ \onen \abs{\hat\tr(B) - \tr(B)} \norm{\beta^{i*}}_{\infty} \\
& \le &
\label{eq::oracle}
D_0 \sigma_{V_i} r_{m,n} + D_0' \tau_B^{1/2}  r_{m,n} (\twonorm{\beta^{i*}} +1)
+ 2 D_1 \inv{\sqrt{m}} \norm{\beta^{i*}}_{\infty} r_{m,n}
\een
By~\eqref{eq::oracle} and the fact that
\bens
2 D_1 := 2\left(\frac{\fnorm{A} }{\sqrt{m}} +\frac{\fnorm{B}  }{\sqrt{n}}
\right)\le  2(\twonorm{A}^{1/2} + \twonorm{B}^{1/2}) (\sqrt{\tau_A} +
\sqrt{\tau_B}) \le D_{\ora} D_0',
\eens
we have  on $\B_0$ and under (A1), for all $i$
\bens
\label{eq::oracleII}
\norm{\hat\gamma^{(i)} - \hat\Gamma^{(i)} \beta^{i*}}_{\infty}
& \le &
D_0' \left( \tau_B^{1/2}  +  \frac{D_{\ora}}{\sqrt{m}} \right) (1+ \twonorm{\beta^{i*}}) r_{m,n}
+ D_0  \sigma_{V_i} r_{m,n} \\
&\le &
C_0 {D_0' \left(\tau_B^{+/2}  \twonorm{\beta^{i*}}
+   \sigma_{V_i} \right) r_{m,n}.}
\eens
Hence the lemma holds for $m \ge 16$ and
$\psi_i = C_0 D_0' K^2 \left(\tau_B^{+/2} \twonorm{\beta^{i*}} + \sigma_{V_i} \right)$.
Finally, we have by the union bound, $\prob{\B_0} \ge 1 -  16/m^2$,
upon adjusting the constant $C_0$.
\end{proofof}

\begin{proofof}{Lemma~\ref{lemma::AD}}
 Lemma~\ref{lemma::AD} follows from the proof of Theorem 26
of~\cite{RZ15}. Hence we only provide a proof sketch.
Recall the following for $X_0 = Z_1 A^{1/2}$,
\bens
\lefteqn{
\Delta := \hat\Gamma -A := \onen X^TX - \onen \hat\tr(B) I_{m} -A} \\
& = & (\onen X_0^T X_0 -A)+  \onen \big(W^T X_0 + X_0^T W\big) + \onen \big(W^T W  - \hat\tr(B) I_{m}\big).
\eens
First notice that
\bens
\maxnorm{\hat\Gamma_A -A}
& \le &
\maxnorm{\onen X_0^T X_0 - A} +
\maxnorm{\onen (W^T X_0 + X_0^T W)} +
\maxnorm{\onen W^T W - \frac{\hat\tr(B)}{n} I_{m}}
\eens
where the last two terms are bounded as in the proof of
Lemma~\ref{lemma::low-noise}.
Putting all together, following the proof of  
Lemma~\ref{lemma::low-noise},
we have~\eqref{eq::defineAD} holds under $\B_{10} \cap \B_{3} \cap \B_{6}$, 
by \eqref{eq::B10a} and \eqref{eq::B10b}, and the fact that 
\bens
\text{ on event} \; \B_3, \quad
\onen \maxnorm{X_0^T X_0 - A} \le 4 C  a_{\max} K^2 \sqrt{{\log 
    m}/{n}},
\eens
where $\prob{\B_{10} \cap \B_{3} \cap \B_{6} }\ge 1 - 6 /m^2 - 6/m^3$;
See Theorem 26 and Corollary 42 of~\cite{RZ15}.
\end{proofof}

\section{Conclusion}
\label{sec::conclude}
The Kronecker sum provides a non-separable covariance model, as an alternative to 
the separable covariance models such as the Kronecker product and may fit 
data better when measurement errors are present.
In this paper, we present theory and method for estimating $\Theta =
A^{-1}$ in the Kronecker Sum covariance model~\eqref{eq::addmodel}. 
Similarly, we can construct the estimator for precision matrix $\Phi = 
B^{-1}$ as the problems for estimating covariance $B$ and $A$ are 
almost symmetrical, 
except that we need to make suitable assumptions to ensure convergence
of both inverse covariance estimators.
Therefore, for both model~\eqref{eq::addmodel}
and~\eqref{eq::matrix-normal-rep-intro}, 
covariance estimation can be obtained through procedures which involve calculating variances for the residual errors after obtaining 
regression coefficients or through MLE refit procedure based on the associated edge set 
\citep[e.g.][]{Yuan10,ZRXB11,Zhou24}.
This procedure can have a broad impact as the measurement errors are 
prevalent in scientific studies. Our methods can be generalized to
numerous scientific data that inherently have measurement errors and
require errors-in-variables models in the high-dimensional setting.
In some sense, we have considered  parsimonious models for 
fitting observation data with two-way dependencies given a single data 
matrix  which is common in genetics and geo-statistical settings. 
The theory developed in this paper can be easily applicable to other
statistical methods that involve the second-order statistics of the
data, i.e. the covariance matrix.

%% file: append.tex
\section*{Organization}
This supplementary material is organized as follows.
In Section \ref{sec:add:proc}, we introduce  an additional estimation
procedure of $\Theta$ as described in Section 2 of the main paper that
guarantees a positive-semidefinite estimator.
Section \ref{Sec:gau} includes details of theoretical properties of
the corrected Lasso estimator from~\cite{RZ15}.
Section \ref{sec:replicate} contains the brief estimation procedure when there are replicates of the data. 
Section \ref{compsec} presents the proposed composite gradient descent
algorithm.
Section \ref{sec:add:real} contains an additional notation table.

\section{Additional estimation procedure}
\label{sec:add:proc}
The estimates $\hat{\Theta}$ and $\hat{\Phi}$ obtained from Algorithm
1 in the main paper are not necessarily positive-semidefinite.
1 \& 2  are not necessarily positive-semidefinite.
To obtain positive-semidefinite estimates, one can consider the following additional estimation procedure.\\

\begin{mdframed}
\noindent
\textbf{Algorithm 3: Obtain $\hat{\Theta}_+$ with an input $\hat{\Theta}$}\\
Consider the case in which $\hat\Theta$ is not positive-semidefinite.
Since $\hat\Theta$ is symmetric,
there exists an orthogonal matrix $U$ and a diagonal matrix $D = \mathrm{diag}(\lambda_1, \cdots, \lambda_m)$ such that
$\hat\Theta = U D U^T$, where $\lambda_1 \le \cdots \le \lambda_m$.
Since $\hat\Theta$ is not positive-semidefinite, $\lambda_1 < 0$.
Let
$
\hat\Theta_+ = U D_+ U^T,
$
where $D_+ :=  \mathrm{diag}(\lambda_1 \vee \epsilon, \cdots, \lambda_m \vee \epsilon)$ for some positive constant $\epsilon \le -\lambda_1$.
Then $\hat\Theta_+$ is positive-semidefinite and satisfies with high probability   
\begin{eqnarray*}
\|\hat\Theta_+ - \Theta\|_2 
\le  \|\hat\Theta - \hat\Theta_+\|_2 + \|\hat\Theta-\Theta\|_2
\le  -2\lambda_1 + \|\hat\Theta-\Theta\|_2
\le 3\|\hat\Theta-\Theta\|_2.
\end{eqnarray*}
\end{mdframed}
The estimate $\hat\Theta_+$ is positive-semidefinite and has the same big-O error bound with $\hat\Theta$ as in Theorem 1 in the main paper.
In practice, we set $\epsilon=-\lambda_1$.

\section{Theoretical properties on EIV regression}
\label{Sec:gau}
Next, we  review theoretical properties of the EIV regression estimator
in Theorem~\ref{thm::lasso}, which is directly from \cite{RZ15}.
We also review the errors-in-variables (EIV) regression model with dependent measurements
studied in \cite{RZ15}.
Suppose that we observe $y\in \R^n$ and $X\in \R^{n\times m}$ in the following regression model
 \begin{equation}
     y = X_0 \beta^* + \epsilon, \qquad \text{where} \qquad X = X_0 + W
 \end{equation}
 where $\beta^* \in \R^m$ is an unknown vector to be estimated, $X_0$
 is an $n\times m$ design matrix and $W$ is a mean zero $n \times m$
 random noise matrix, independent of $X_0$ and $\epsilon$, whose
 columns are also independent and each consists of dependent
 elements.
 That is, we consider $\E \omega^j \otimes \omega^j = B$ for 
 all $j=1, \dots, m$, where $\omega^j$ denotes the $j^{th}$ column 
 vector of $W$; cf. Definition~\ref{def::subgdata}.
  We assume that the noise vector $\ve \in \R^n$ is independent of $W$
  or $X_0$, with independent entries $\ve_{j}$ satisfying $\E[\ve_{j}]
  = 0$ and $\norm{\ve_{j}}_{\psi_2} \leq M_{\ve}$, where
  $\norm{x}_{\psi_2} := \mathrm{inf}\{s>0 \mid \mathbb{E}
  [e^{(x/s)^2}] \le 2\}$ for  a scalar random variable $x$.
A random variable $Z$ is sub-gaussian if and only if $S :=
Z^2$ is sub-exponential with $\norm{S}_{\psi_1} =
\norm{Z}^2_{\psi_2}$. More generally, we have 
Lemma~\ref{lemma::product}.

\begin{definition}
\label{def::subgdata}
Let $Z$ be an $n \times m$ random matrix with independent entries $Z_{ij}$ satisfying
$\E Z_{ij} = 0$, $1 = \E Z_{ij}^2 \le \norm{Z_{ij}}_{\psi_2} \leq K$.
Let $Z_1, Z_2$ be independent copies of $Z$.
Let $X = X_0 + W$ such that
\bnum
\item
$X_0 = Z_1 A^{1/2}$ is the design matrix with independent subgaussian
row vectors,
\item
$W =B^{1/2} Z_2$ is a random noise matrix with independent subgaussian
column vectors.
\enum
\end{definition}

Before we continue, we need to 
define some parameters related to the restricted and
sparse eigenvalue conditions that are needed to state our main results
on EIV regression.
We first state Definitions~\ref{def:memory}-~\ref{def::lowRE}.
For more details of these, see \citet{RZ15}.
\begin{definition}
\label{def:memory}
\textnormal{\bf (Restricted eigenvalue condition $\RE(s_0, k_0, A)$)}
Let $1 \leq s_0 \leq m$ and  $k_0$ be a positive number.
The $m \times m$ matrix $A$ satisfies $\RE(s_0, k_0, A)$
 condition with parameter $K(s_0, k_0, A)$ if for any $\upsilon
 \not=0$,
\beq
\inv{K(s_0, k_0, A)} :=
\min_{\stackrel{J \subseteq \{1, \ldots, p\},}{|J| \leq s_0}}
\min_{\norm{\upsilon_{J^c}}_1 \leq k_0 \norm{\upsilon_{J}}_1}
\; \;  \frac{\norm{A \upsilon}_2}{\norm{\upsilon_{J}}_2} > 0,
\eeq
where $\upsilon_{J}$ represents the subvector of $\upsilon \in \R^m$
confined to a subset $J$ of $\{1, \ldots, m\}$.
\end{definition}

\begin{definition}{\textnormal{(Upper-$\RE$ condition)}}
\label{def::upRE}
The $m \times m$ matrix  $\Gamma$ satisfies an upper-$\RE$ condition with curvature
$\bar\alpha >0$ and tolerance $\tau > 0$ if
\bens
\theta^T \Gamma \theta \le \bar\alpha \twonorm{\theta}^2 + \tau
\|\theta\|_1^2, \; \;  \forall \theta \in \R^m.
\eens
\end{definition}

\begin{definition}
\label{def::sparse-eigen}
Define the largest and smallest
$d$-sparse eigenvalue of a $m \times m$ matrix $A$: for $d<m$,
\ben
\label{eq::eigen-Sigma}
\rho_{\max}(d, A)  :=
\max_{t \not= 0; \|t\|_0 \le d} \; \;\frac{\shtwonorm{A
    t}^2}{\twonorm{t}^2}, \quad
\rho_{\min}(d, A)  :=
\min_{t \not= 0; \|t\|_0 \le d} \; \;\frac{\shtwonorm{A t}^2}{\twonorm{t}^2}.
\een
\end{definition}

Throughout this section, let $s_0 \ge 1$ be the largest integer chosen such that the following inequality holds:
\ben
 \label{eq::s0cond}
\sqrt{s_0} \vp(s_0) \le \frac{\lambda _{\min}(A)}{32 C}\sqrt{\frac{n }{\log m}}, \quad
\; \; \vp(s_0) := \rho_{\max}(s_0, A)+\tau_B,
\een
where  $\tau_B = \tr(B)/n$ and $C$ is to be defined.
{
Denote by
\ben
\label{eq::defineM}
M_A = \frac{64 C \vp(s_0)}{\lambda_{\min}(A)} \ge 64 C.
\een}
We use the expression
$\tau := (\lambda_{\min}(A)-\alpha)/{s_0}, \; \; \text{where} \; \; \alpha =
5\lambda_{\min}(A)/8.$

\subsection{EIV regression}



\begin{theorem}{\textnormal{(\bf{Estimation for the Lasso-type estimator})}}[Theorem 3 of \citet{RZ15}]
\label{thm::lasso}
Suppose $n=\Omega(\log m)$ and $n \le (\V/e) m \log m$,
where $\V$ is a constant which depends on $\lambda_{\min}(A)$,
$\rho_{\max}(s_0, A)$ and $\tr(B)/n$.
Suppose $m$ is sufficiently large. Suppose (A1)-(A4) hold.
Consider the EIV regression model (1) and (3) as defined in the main paper
with independent random matrices $X_0, W$ as in
Definition~\ref{def::subgdata} with $\norm{\e_{j}}_{\psi_2} \leq M_{\e}$.
Let $C_0, c' > 0$ be  some absolute constants. Let $D_2 := 2(\twonorm{A} + \twonorm{B})$.
Suppose that $c' K^4 \le 1$ and
\ben
\label{eq::trBLasso}
\quad
r(B) := \frac{\tr(B)}{\twonorm{B}} & \ge & 16c' K^4 \frac{n}{\log m}
\log \frac{\V m \log m }{n}.
\een
Let $b_0, \phi$ be numbers which satisfy
\ben
\label{eq::snrcond}
\frac{M^2_{\e}}{K^2 b_0^2}   \le \phi  \le 1.
\een
Assume that the sparsity of $\beta^*$ satisfies for some $0 < \phi \le
1$
\ben
\label{eq::dlasso}
&& d:= \abs{\supp(\beta^*)} \le
\frac{c' \phi K^4}{40 M_+^2} \frac{n}{\log  m}< n/2,\\
&&\mathrm{where}\  M_+ = \frac{32C\varpi(s_0+1)}{\lambda_{\min}(A)}
\een
for $\varpi(s_0+1) = \rho_{\max}(s_0+1,A) + \tau_B$.
Let $\hat\beta$ be an optimal solution to the EIV regression  with
\ben
\label{eq::psijune}
&& \lambda \ge 4 \psi \sqrt{\frac{\log m}{n}} \; \; \text{ where } \;\;
\psi  := C_0 D_2 K \left(K \twonorm{\beta^*}+ M_{\e}\right).
\een
Then for any $d$-sparse vectors $\beta^* \in \R^m$, such that
$\phi b_0^2 \le \twonorm{\beta^*}^2 \le b_0^2$, we have with probability  at least $1- 16/m^3$,
\bens
\twonorm{\hat{\beta} -\beta^*} \leq \frac{20}{\alpha}  \lambda \sqrt{d} \; \;
\text{ and } \; \norm{\hat{\beta} -\beta^*}_1 \leq \frac{80}{\alpha}
\lambda d.
\eens
\end{theorem}

\subsection{Multiple EIV regressions}

The following theorem \ref{thm::nodewise} shows oracle inequalities of the nodewise regressions, which is analogous to Theorem 6 of \citet{RZ15}.
\begin{theorem}
\label{thm::nodewise}
Consider the Kronecker sum model as in \eqref{eq::addmodel}.
Suppose all conditions in Theorem~\ref{thm::lasso} hold, except that
we drop \eqref{eq::snrcond} and replace~\eqref{eq::psijune} with
\begin{equation}
\label{eq::psijune15}
 \lambda^{(i)} \ge 4 \psi_i\sqrt{\frac{\log m}{n}} \; \; \text{ where } \;\;
\psi_i := C_0 D_0' K^2  \left(\tau_B^{+/2} \twonorm{\beta^{i*}} + \sigma_{V_i}  \right),
\end{equation}
where $\sigma_{V_i}^2 := A_{ii}-A_{i,-i} A_{-i,-i}^{-1} A_{-i,i}$.    
Suppose that for $0< \phi \le 1$ and $C_A := \inv{160 M_{+}^2}$,
\ben
\label{eq::doracle}
&& d:= \max_{1 \le i \le m}\abs{\supp(\beta^{i*})} \le C_A \frac{n}{\log m} \left\{c' c'' D_{\phi}
  \wedge 8 \right\},  \;\;\text{ where }  \\
\nonumber
&&    
c'' := \frac{\twonorm{B} + a_{\max}}{\varpi(s_0+1)^2} , \;
 D_{\phi}
 = \frac{K^2 M^2_{\e}}{b_0^2}  + K^4 \tau_B^{+}  \phi,   
\een
and $c', \phi, b_0$ 
as defined in
Theorem~\ref{thm::lasso}.
We obtain $m$ vectors of  $\hat{\beta}^i, i=1, \ldots, m$ 
by solving (7) in the main paper with $\lambda = \lambda^{(i)}$,
\bens
 \; \hat\Gamma^{(i)} & = &\onen X_{\minus i}^T  X_{\minus i} - \hat\tau_B I_{m-1}
\; \text{ and } \; \hat\gamma^{(i)}
 \; = \; \onen  X_{\minus i}^T X_i,
\eens
for each $i$. 
Let $b_1 := b_0 \sqrt{d}$.
Then for all $d$-sparse vectors  $\beta^{i*} \in \R^{m-1}, i=1, \ldots, m$, such that
$\phi b_0^2 \le \twonorm{\beta^{i*}}^2 \le b_0^2$,
we have with probability  at least $1- 32/m^2$,
\bens
\twonorm{\hat{\beta}^{i} -\beta^{i*}} \leq \frac{20}{\alpha}
\lambda^{(i)}  \sqrt{d} \; \;
\text{ and } \; \norm{\hat{\beta}^{i} -\beta^{i*}}_1 \leq \frac{80}{\alpha}
\lambda^{(i)} d.
\eens
\end{theorem}


Theorem \ref{coro::Omega} shows the 
consistency of $\hat{\Phi}$ in the operator norm.

\begin{theorem}
\label{coro::Omega}
Suppose the columns of $\Phi$ is $d_{\Phi}$-sparse, and suppose the
condition number $\kappa(\Phi)$ is finite.
Suppose conditions of Theorem 3 in the main paper hold upon
swapping $m$ for $n$ and vice versa, and replacing $d$ with $d_{\Phi}$.
Then with probability at least $1- 26 /m^2$,
\bens
\twonorm{\hat\Phi - \Phi}  =  O_p\left( \frac{K^2(B)}{\lambda_{\min}^2(B)}  d_{\Phi} \max_{i \le n} \lambda^{(i)} \right).
\eens
\end{theorem}

 \section{Using replicates}
\label{sec:replicate}
In the current work, we focus on a class of generative models which rely
on the sum of Kronecker product covariance matrices to model complex trial-wise
dependencies as well as to provide a general statistical framework for
dealing with signal and noise decomposition. 
It can be challenging to handle this type of data because 
there are often 
dependencies among the trials, and
thus they can not be treated as independent replicates.
In the present work, we propose a framework
for explicitly modeling the variation in a set of ``replicates''. 
We denote the number of subjects by $n$,  the number of time
points by $m$, and the number of  replicates by $N$.  
Let $A \in \R^{m \times m}$ and $B, C  \in \R^{n\times n}$ be $m$ by $m$ and $n$ by $n$ positive definite 
matrices, respectively. We consider the following generative model:
\begin{eqnarray}
\label{main:eq:1}
&& X_i= X_{0,i}+W_i, \quad \forall i=1,\cdots, N 
\; \; \text{where }\\
\nonumber
&& \mathrm{vec}(X_{0,i}) \sim \mathcal{L}(0, A \otimes C) \quad \text{ and
}  
\quad \mathrm{vec}(W_i) \sim \mathcal{L}(0, I_m \otimes B),
\end{eqnarray}
which means that $X_{0,i}, W_i \in \R^{n \times m}$ are independent subgaussian random
matrices such that $\mvec{X_{0,i}}, \mvec{W_i} \in \R^{mn}$ 
have covariances $A \otimes C$ and $I_m \otimes B$, respectively.
Here the mean response matrix $X_{0,i}$ and the experiment-specific variation matrices $W_i$
can jointly encode the temporal and spatial dependencies. 
More generally, we can use a subgaussian random matrix
to model replicate-to-replicate fluctuations:
\bens
\mvec{W_i} \sim \mathcal{L}(0, I_m \otimes B(t)),\; \; \text{ where } \; \; B(t) \succ 0
\eens
is the covariance matrix describing the spatial dependencies which may
vary cross the $N$ replicated experiments.
In this subsection, we consider the general model \eqref{main:eq:1},
but with 
\begin{equation}
\label{eq:main}
\mathrm{vec}(X_{0,i}) \sim \N(0, A \otimes I_n) \; \text{ and } \; 
 \mathrm{vec}(W_i) \sim \N(0, I_m \otimes B),
\end{equation}
where $X_{0,i}$ and $W_i$ are independent of each other.

\noindent
\textbf{Case 1}: if $X_{0,i} = X_0$,  one can avoid the assumption that
the trace of $A$ is known.
To estimate $\Phi =B^{-1}$, we note that
\[
\forall i \neq j,\quad
\mathrm{vec}(X_i-X_j) \sim \N(0, I_m \otimes 2B).
\]
Without loss of generality, assume $N$ is an even number. We have $Nm/2$ replicates to estimate $\Phi$: each of the $m$ columns of $\widetilde{W}_i=X_{2i-1}-X_{2i}$ for $i=1,\cdots, N/2$ is a random sample of $\N(0_n, 2B)$.    
Then by using nodewise regression method \citep{Yuan10,ZRXB11} as described in Section \ref{sec:graphical}, we can estimate $\Phi$ with the Frobenius error rate of
$O_p \left( d_{\Phi} \sqrt{\frac{\log n}{Nm}}\right)$, where $d_{\Phi}$ is the sparsity parameter of each column of $\Phi$, i.e., 
there exists at most $d_{\Phi}$ non-zero components in each column of $\Phi$.

To estimate $\Theta=A^{-1}$, consider the following observed mean response:
\ben
\label{eq::Beffect}
\overline{X} =\frac{1}{N} \sum_{i=1}^N X_i, \ \rm{where}\
\cov(\mvec{\overline{X}}) =  A \otimes I_n + \frac{1}{N} I_m \otimes B.
\een
Note that when $N$ is sufficiently large, then $\overline{X}$ is close to the traditional covariance structure having independent rows.
By using $\hat\tau_B  := \onen \hat\tr(\hat{B})$,  one can show $|\hat\tau_B- \tau_B| = O_p \left( d_{\Phi} \sqrt{\frac{\log n}{Nm}}\right)$, then
we obtain the Frobenius error bound of $\hat{\Theta}$ as in Theorem \ref{coro::Theta}.
We omit detailed proofs as they utilize existing results.  We
summarize the estimation procedures as follows: \\

\noindent
\textbf{Algorithm 2-1: Obtain $\hat{\Phi}$ with i.i.d. input vectors $\{\widetilde{W}_i\}_{i=1,\cdots, N/2}$}\\
\noindent
\noindent{\textbf{(1)}}
Let $\widetilde{B} := \sum_{i=1}^{N/2} \widetilde{W}_i
\widetilde{W}_i^T /(Nm) \succeq 0$ be an unbiased estimator of $B$.\\
\noindent{\textbf{(2)}}
Apply graphical Lasso \citep{FHT07} with the input $\widetilde{B} \succeq 0$ to obtain $\hat{\Phi}$.\\
\noindent
\textbf{Algorithm 2-2: Obtain $\hat{\Theta}$ with input $\hat{\tau}_B=
  \tr(\widetilde{B})/n$}.\\  
\noindent
Estimate $\Theta$ using  Algorithm 1 with $$\hat \Gamma = \overline{X}^T \overline{X}/n-\frac{\hat{\tau}_B}{N} I_m,\quad
\hat\Gamma^{(j)} =\hat \Gamma_{-j, -j} 
\quad \text{ and } \; \hat\gamma^{(j)} \; = \; \onen  \overline{X}_{\minus j}^T \overline{X}_j.
$$
 
\noindent
\textbf{Case 2}: if $X_{0,i}$ are independent, note that 
\begin{align*}
& \tilde X = [X_1,\ldots, X_N] \in \R^{m \times Nn}, \ \rm{where}\
\cov(\mvec{\tilde{X}}) =   (I_N \otimes A) \oplus B \\
& \ddot X = [X_1^\top,\ldots, X_N^\top]^\top \in \R^{Nm \times n}, \ \rm{where}\
\cov(\mvec{\ddot{X}}) =   A \oplus (I_N \otimes B). 
\end{align*}
We can estimate $\Theta = A^{-1}$ and $\Omega = B^{-1}$ from $\ddot X$ and $\tilde X$ by applying the proposed Algorithm 1, respectively.

\vskip 1cm

\section{Composite gradient descent algorithm}
\label{compsec}
For a function $g: \R^m \to \R$, we write $\grad g$ to denote a 
gradient or subgradient, if it exists.
To solve the non-convex optimization problem~\eqref{eq::origin} in Algorithm 1,   
we use the composite gradient descent algorithm as studied in \citet{ANW12} and \citet{LW12}. Let 
$L(\beta):=    \frac{1}{2} \beta^T \hat\Gamma^{(i)} \beta - \ip{\hat\gamma^{(i)}, \beta}$.
The gradient of the loss function is $\grad L(\beta) = \hat{\Gamma}^{(i)} \beta -\hat\gamma^{(i)}$. The composite gradient
descent algorithm produces a sequence of iterates $\{\beta^{(t)},\ t=0,1,2, \cdots, \}$ by
\begin{equation}
\label{grad_comp_2}
\beta^{(t+1)}= \argmin_{\|\beta\|_1 \le b_1} L(\beta^{(t)})+ \langle \grad L(\beta^{(t)}), \beta- \beta^{(t)} \rangle + \frac{\eta}{2} \|\beta-\beta^{(t)}\|_2^2 + \lambda \|\beta\|_1 
\end{equation}
with the step size parameter $\eta>0$.

\section{Additional Table}
\label{sec:add:real}
In this section we include the table of notations used throughout the manuscript.

\begin{table}[H]
\caption{\label{table_def}The Symbols}%
\centering
\fbox{
\begin{tabular}{l  c } 
Parameters   & Definitions \\  [0.5ex] 
\hline 
 $X_{i}$ &  The $i$th column of a matrix $X$ \\
 $X_{\minus i}$ &  The sub-matrix of $X$ without the $i$th column \\
 $X_{\minus i, \minus j}$ &  The sub-matrix of $X$ without the $i$th row and $j$th column \\
$\tau_{\Sigma}$ &   $\rm{tr}(\Sigma)/p$ for a square matrix $\Sigma \in \R^{p \times p}$\\
$\Theta$ & $A^{-1}$\\
$\Phi$ & $B^{-1}$\\
$a_{\max}$ & $\max_i A_{ii}$ \\
$b_{\max}$ & $\max_i B_{ii}$ \\
$\alpha$ & $\frac{5}{8} \lambda_{\min}(A)$\\
$\bar{\alpha}$ & $\bar{\alpha}=8 \lambda_{\max}(A)/11$\\
$\tau$ & $(\lambda_{\min}(A)-\alpha)/{s_0}$\\
$\eta \ge \frac{11}{8}\lambda_{\max}(A)$ & step size parameter \\
$\vp(s_0)$  &   $\rho_{\max}(s_0, A)+\tau_B$ \\
$\tau_0$ &  $\frac{400C^2 \vp(s_0+1)^2}{\lambda_{\min}(A)}$\\
$D'_0$ &  ${\twonorm{B}}^{1/2} + a_{\max}^{1/2}$\\
$\tilde{D}'_0$  & ${\twonorm{A}}^{1/2} + b_{\max}^{1/2}$\\
$D_{\ora}$  &   $2(\twonorm{A}^{1/2} + \twonorm{B}^{1/2})$\\
$\tau_B^{+/2}$  & $\sqrt{\tau_B} + \frac{D_{\ora}}{\sqrt{m}}$\\
$C_\psi$ & $ K^2 C_0 D_0'  \left(\tau_B^{+/2} \kappa(A)+ \sqrt{a_{\max}} \right)$ \\
$\psi_0$ &   $0.1 D_0' \left(\tau_B^{+/2} a_{\max}+ \sqrt{a_{\max}}\right)$\\
$\rho_{\max}(d, A)$ &  $\max_{t \not= 0; d-\text{sparse}} \; \;\shtwonorm{A
  t}^2/\twonorm{t}^2$\\
$\vp(s_0)$ & $\rho_{\max}(s_0, A)+\tau_B$\\
$s_0$  &    The largest integer satisfying $\sqrt{s_0} \vp(s_0) \le \frac{\lambda _{\min}(A)}{32 C}\sqrt{\frac{n }{\log m}}$\\
$\psi_1^2$ &   $\frac{16c}{\lambda_{\min}(A)}\left(\frac{1}{1-\kappa}+\frac{\lambda_{\min}(A)}{2s_0}  \right)$\\
\end{tabular}
}
\end{table}